\documentclass{article}
\usepackage{arxiv}
\usepackage{graphicx,amsmath,amsfonts,amscd,amsthm,amssymb,pstricks}
\usepackage{color}
\usepackage{dcolumn}
\usepackage{bm}
\usepackage{longtable}
\usepackage{mathrsfs}
\usepackage{textcomp}
\usepackage{esvect}
\usepackage{float}
\usepackage{mathtools,amsfonts}
\usepackage[USenglish,british]{babel}
\usepackage{environ}
\usepackage{xifthen}
\usepackage{xargs}		
\usepackage{hyperref}
\usepackage{physics,mathtools}
\usepackage[separate-uncertainty = true]{siunitx}
\usepackage{multirow}
\usepackage{comment}
\usepackage{enumerate}
\usepackage{appendix}
\usepackage{tikz}	
\usetikzlibrary{arrows, calc, matrix, patterns, decorations.markings, shapes, decorations.pathmorphing, shadows.blur}
\usepackage[utf8]{inputenc}
\usepackage{enumitem}

\newcommand{\av}[1]{\left<#1\right>}
\newcommand{\ie}{i.e.\ }

\newcommand{\eg}{e.g.\ }

\newcommand{\wrt}{w.r.t.\ }

\newcommand{\ham}{\mathcal{H}}

\newcommand{\eps}{\varepsilon}

\newcommand{\omegar}{\omega_\mathrm{r}}

\title{Analysis of double-resonance crossing in adiabatic trapping phenomena for quasi-integrable area-preserving maps with time-dependent exciters}

\author{A. Bazzani\\
Physics and Astronomy Department, Bologna University and INFN-Bologna, V. Irnerio 46, 40126 Bologna - Italy\\
\And
F. Capoani\\
Physics and Astronomy Department, Bologna University and INFN-Bologna, V. Irnerio 46, 40126 Bologna - Italy\\
\And
M. Giovannozzi\thanks{Corresponding author: massimo.giovannozzi@cern.ch}\\
Beams Department, CERN, Esplanade\ des Particules 1, 1211 Geneva 23, Switzerland}

\begin{document}
\maketitle

\begin{abstract}
In this paper, we analyze the adiabatic crossing of a resonance for Hamiltonian systems when a double-resonance condition is satisfied by the linear frequency at an elliptic fixed point. We discuss in detail the phase-space structure on a class of Hamiltonians and area-preserving maps with an elliptic fixed point in the presence of a time-dependent exciter. Various regimes have been identified and carefully studied. This study extends results obtained recently for the trapping and transport phenomena for periodically perturbed Hamiltonian systems, and it could have relevant applications in the adiabatic beam splitting in accelerator physics. 
\end{abstract}



\section{Introduction}

The adiabatic theory for Hamiltonian systems addresses the problem of understanding the consequences of slow parametric modulations. The concept of adiabatic invariance allows one to predict the long-term evolution of a system, highlighting the fundamental properties of action variables~\cite{arnold,chirikov}. The theory has been fully developed in the case of one degree-of-freedom systems~\cite{neish1976,neish1985,cary1986,cary1988,elskens1991,elskens1993,neish1997,an6,dv_im,an10,neish2008,mosovsky2011,neish2021}, and its extension to the multidimensional case or to symplectic maps~\cite{ab_fb_gt_AIP} is a difficult problem, mainly due to the small denominators of perturbation theory and the ubiquitous presence of nonlinear resonances in phase space~\cite{via1,via3}. 

The possibility of manipulating the phase-space structure in an adiabatic way has recently been considered for novel applications in the realms of accelerator and plasma physics~\cite{an4,sridhart,bvc4,escande2016,mte1,mte2}. For instance, nonlinear resonance trapping and adiabatic transport have been employed to manipulate a charged-particle distribution to minimize particle losses during the extraction process of a circular accelerator of high-intensity beams. Adiabatic manipulations are also used to provide control of transverse beam emittances~\cite{mte2017,mte3,mte4}. These experimental procedures~\cite{mte2017,mte3,mte4} require very precise control of the efficiency of adiabatic trapping in resonances~\cite{an10,an9,add2}, as well as the change in phase space during adiabatic transport. All these processes can be modeled using multidimensional Hamiltonians or symplectic maps~\cite{bazzanietal,Bazzani2022}. 

An interesting and intriguing observation has been made during the experiments carried out at the CERN Proton Synchrotron (PS) for Multi-Turn Extraction (MTE)~\cite{Borburgh:2137954,Abernethy2017,mte2017,Huschauer2017}. Experimental observations have clearly indicated that the efficiency of the beam trapping into a nonlinear resonance can be improved by using an external exciter whose frequency is set on resonance condition with the main frequency of the system. The resulting model is a quasi-resonant Hamiltonian system perturbed by a time-dependent external exciter whose frequency satisfies a double-resonance condition. We have studied in detail the phase-space structure of this system considering the adiabatic crossing of the resonance in different dynamical regimes. The results presented in the article extend the recent results presented in Ref.~\cite{Bazzani2022}. 

The paper is organized as follows: In Section~\ref{sec:recap} the main results of the theory of adiabatic invariance and separatrix crossing are briefly summarized. In Section~\ref{sec:models} we discuss two models: A Hamiltonian model is introduced and used for carrying out the analytical computations and to understand the dependence of the phase-space structure on the system's parameters; a map model is used for numerical simulations because it represents a more realistic model of dynamics. The Hamiltonian system allows one to identify the phase-space topologies and the various regimes of the system, and numerical simulations verify the robustness of the analytical results. The analysis of the trapping process is discussed in Section~\ref{sec:trapping}, while numerical simulations are presented and discussed in detail in Section~\ref{sec:simresults}. Finally, some conclusions are drawn in Section~\ref{sec:conc}, while some detailed computations can be found in the appendices. 
\section{Adiabatic theory for trapping in a stable resonance}\label{sec:recap}

The results of this paper take advantage of the theory of separatrix crossing in adiabatic conditions, which describes how an orbit may adiabatically cross a separatrix, breaking the adiabatic invariance of the action in a controlled way. Here, we recall the main results that will be useful for the analysis carried out in this paper. 

According to Neishtadt's theory of adiabatic trapping~\cite{neish1975}, which can also be applied to area-preserving map models~\cite{bazzanietal,Bazzani2022}, when a slow modulation of a parameter makes the area of a separatrix-enclosed region of phase space equal to $2\pi J_0$, where $J_0$ is the initial action variable, the particle jumps into another region of the phase space with a probability that depends on the time derivative of the phase-space regions.

\begin{figure}[htb]
  \begin{center}
    \includegraphics[trim=45truemm 185truemm 45truemm 45truemm,width=0.5\columnwidth,clip=]{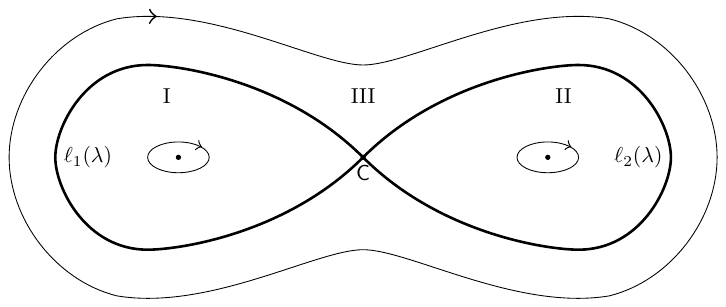}
  \end{center}
  \caption{Generic phase-space portrait divided into three regions ($\mathrm{I}$, $\mathrm{II}$, $\mathrm{III}$) by the separatrices $\ell_1(\lambda)$ and $\ell_2(\lambda)$.}
  \label{fig:genericphsp}
\end{figure}

In general, consider a Hamiltonian system dependent on a slowly modulated parameter $\lambda=\epsilon t, \epsilon \ll 1$, whose phase space is divided by two separatrices $\ell_1(\lambda)$ and $\ell_2(\lambda)$, into three regions denoted as Region~$\mathrm{I}$, $\mathrm{II}$ and $\mathrm{III}$ according to the sketch in Fig.~\ref{fig:genericphsp}. Let us indicate by $A_i(\lambda)$ the area of the region $i$ for a certain value of $\lambda$.

If an orbit of an initial condition in Region~$\mathrm{III}$ encloses an area $A_0=2\pi J_0$, where the action $J_0$ is the adiabatic invariant, when $\lambda$ assumes a value so that $A_\text{III}(\tilde\lambda) = A_0$, then the particle will enter Region $\text{I}$ or $\text{II}$ according to the probability $\mathcal{P}_{\mathrm{III}\to i}$ given by the formula:
\begin{equation}
\mathcal{P}_{\mathrm{III}\to i} =
\begin{dcases}
0\qquad & \text{if } \xi_i < 0\\
\xi_i\qquad & \text{if } 0\le \xi_i \le 1\\
1\qquad & \text{if }  \xi_i > 1
\end{dcases}\qquad i=1,2 \, , 
\end{equation}
where
\begin{equation} 
\xi_i = \frac{\dv*{{A}_\mathrm{i}}{\lambda}}{\dv*{A_\mathrm{I}}{\lambda}+\dv*{A_\mathrm{II}}{\lambda}}\eval_{\lambda=\tilde\lambda}\, .
 \label{eq:neish}
\end{equation}

\section{The dynamical models} \label{sec:models}
\subsection{The Map model}
We consider a modified version of the Hénon map~\cite{henon} by adding a cubic nonlinearity and a modulated exciter
\begin{equation}
\begin{split}
    \begin{pmatrix}
    x_{n+1} \\ p_{n+1} 
    \end{pmatrix} & = 
    R(\omega_{0,n}) \times \\
    & \times \begin{pmatrix}
    x_n\\
    p_n+x_n^2+\kappa x_n^3 + \eps_n\cos(\omega_{n} n + \psi_0)
    \end{pmatrix} \, , 
\end{split}
    \label{eq:doubleresmap}
\end{equation}
where $R(\omega_{0,n})$ represents a time-dependent rotation matrix. The reason for adding a cubic nonlinearity is twofold: It makes the system closer to that used in the application~\cite{Borburgh:2137954,Abernethy2017,mte2017,Huschauer2017}; it makes it possible to build an interpolating Hamiltonian that, even at the lowest order of perturbation theory, is capable of describing the topological structure of the phase space. The external exciter is customarily described by its frequency $\omega_n$, its strength $\eps$, and its phase $\psi_0$. All three parameters can be set to be time-dependent, but in our studies we only vary the first two. 

\begin{figure*}
\includegraphics[width=\textwidth]{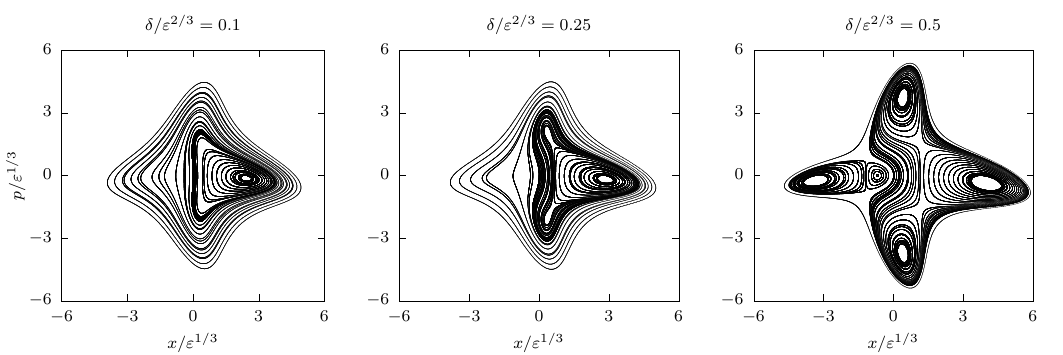}
\caption{Phase-space portraits of the Poincaré map of Eq.~\eqref{eq:doubleresmap} sampled every $4$ iterations. We use three values of $\delta/\eps^{2/3}$ that account for three possible resonance topologies (the other parameters values are: $\omegar/(2\pi)=1/4$, $\Delta=0$, $\eps=10^{-4}$, $\kappa=\num{0.1}$, $\psi_0=0$).}
\label{fig:phasespace_map}
\end{figure*}

Our aim is to study the system under the effect of two simultaneous resonant conditions. Setting $\omegar\in 2\pi\mathbb{Q}$ as the resonant frequency, we choose the two frequencies so that $\omega_{0,n}\approx \omega_n \approx \omegar$, and we introduce the two parameters $\delta_n$ and $\Delta_n$ to account for the distance to each resonance, \ie  $\omega_{0,n} = \omegar + \delta_n$ and $\omega_{n} = \omegar + \Delta_n$. The three parameters $\delta_n$, $\Delta_n$, and $\eps_n$ depend on the iteration number, being slowly modulated to satisfy an adiabatic condition.

We consider the case close to the resonance $1:4$, \ie $\omegar/(2\pi)=1/4$, and, according to the Poincaré-Birkhoff theorem~\cite{arnoldergo}, we expect that a chain of four resonance islands enclosing a central core appears in phase space for certain values of $\delta$. However, the effect of the second resonance condition $1:1$, generated by the external exciter, changes the structure of the phase space. To study these effects, we set $\Delta_n=0$, so that the excitation frequency is exactly resonant with $\omegar$ and analyze the phase-space structure of the frozen map (i.e.,$\omega_{0,n}$ and $\omega_n$ are kept constant) using the $4$th iterate of Eq.~\eqref{eq:doubleresmap}, which is shown in Fig.~\ref{fig:phasespace_map}. Of course, when $\eps_n=0$ one retrieves the well-known four-island structure of the Hénon map~\cite{henon,yellowreport}. The second resonance acts on the island structure by enlarging the size of one or two selected islands and reducing the size of the others. Furthermore, when $\eps$ is large enough, other fixed points disappear, destroying one or more islands. 

\begin{figure}
\centering
\includegraphics[width=0.5\columnwidth]{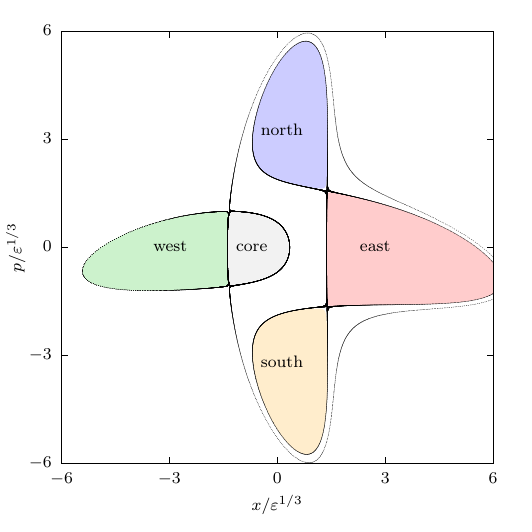}
\caption{Separatrices (black lines) of the Poincaré map of Eq.~\eqref{eq:doubleresmap} sampled every $4$ iterations. The four islands and the core have been filled with the colors used in Figs.~\ref{fig:area_map_eps}, \ref{fig:areadiff}, \ref{fig:plot_distr_map} and \ref{fig:plot_eps_map} to refer to each region. The naming convention of the various regions is also reported here (parameters values: $\omegar/(2\pi)=1/4$, $\Delta=0$, $\eps=10^{-4}$, $\kappa=\num{0.1}$, $\psi_0=0$ and $\delta/\eps^{2/3}=1$).}
\label{fig:phsp_key}
\end{figure}

Following the sketch in Fig.~\ref{fig:phsp_key}, which represents a situation where all four islands and the core are present, we outline the identification of each island by means of cardinal points according to their position in the phase space $(x,p)$. Furthermore, the color used in this figure corresponds to the convention used in subsequent figures, for the sake of clarity. Notable is the peculiar structure of the separatrices: Without the exciter, there exists a single separatrix that connects the hyperbolic fixed points and encloses the elliptic fixed points. In the presence of the exciter, the single separatrix splits into two structures: One divides the core region from the West island and encloses a region of the phase space that contains the other three islands; the other encloses the remaining three islands but does not include the core region. As $\eps \to 0$, the separatrix that surrounds the three islands and the outer separatrix merge together.

We also observe that increasing $\eps$, for $\psi_0=0$, the East island grows, whereas the surface of the other islands shrinks (the island opposite the dominant, the West island in this case, becomes even smaller). Of course, since the exciter frequency is exactly $1/4$, the islands remain fixed in phase space only when we consider the stroboscopic map (i.e.,the $4$th iterate of the original map). If we were to observe each iteration of the map, we would see islands rotating each turn by $\pi/2$, so that the main island, depending on the value of $n\pmod 4$, can be found at each cardinal point.

\begin{figure*}
\includegraphics[width=\textwidth]{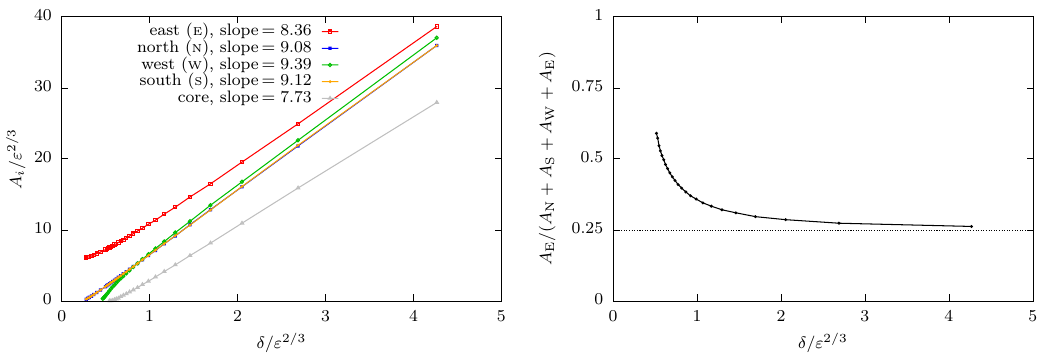}%
\caption{Left: Plot of the area $A_i$ of each island and of the core region in the $4$th iterate of the map of Eq.~\eqref{eq:doubleresmap} as a function of $\delta/\eps^{2/3}$. Note that the lines for the North and South islands are almost perfectly overlying. The slope quoted is the angular coefficient of the linear fit of each line performed in the interval $1\le \delta/\eps^{2/3} \le 4$. Right: Ratio between the area of the main island (East) and the sum of all islands using the data of the left plot (parameters values: $\omegar/(2 \pi)=1/4$, $\Delta=0$, $\kappa=\num{0.1}$, $\psi_0=0$).}
\label{fig:area_map_eps}
\end{figure*}

The measure of the island area as a function of $\delta/\eps^{2/3}$ is shown in Fig.~\ref{fig:area_map_eps} (left), while in the right plot we show the ratio between the area of the main island and the sum of the areas of all islands. We observe that as $\eps\to 0$ or $\delta\to\infty$ one recovers the $1/4$ ratio as expected in the case of the modified Hénon map without the exciter. On the contrary, for small values of $\delta$ and large values of $\eps$, the main island becomes larger than the others. Furthermore, from the plot on the left, one observes that the variation of $A_i$ as a function of $\delta/\eps^{2/3}$ becomes linear for large values of the parameter. Note also that the choice of this combination of model parameters $\delta$ and $\eps$ is justified by the analysis of the Normal Form Hamiltonian that is discussed in Section~\ref{sec:ham}.

The same figure highlights two critical values of the re-scaled parameter. The East, South and North islands are created for a positive value of $\delta/\eps^{2/3}$, which corresponds to a bifurcation point when the separatrix enclosing these islands is created. For lower values of $\delta/\eps^{2/3}$, the phase-space structure that morphs into the East island is the only stable region whose center is displaced from the origin by the exciter's action, and therefore there is no island present. Then, there is a second critical value of $\delta/\eps^{2/3}$ for the other two regions (i.e., the core and the West island), which are defined by the second inner separatrix. These critical values will be discussed in conjunction with the Hamiltonian model in Section~\ref{sec:ham}.

Although it is the largest island, the variation in the surface of the East island is not the largest. As resonance trapping is determined by their area derivative, the majority of particles found in that island after the modulation process comes from the displacement of the center rather than from the separatrix crossing.

If $\eps$ is fixed and $\delta$ varies, all island areas increase, but the East island remains the largest. The area ratio between the East island and the sum of the four islands is higher for small values of $\delta$ and tends to $1/4$ (i.e., the islands enclose the same area) as the frequency moves farther away from the resonance (see right plot of Fig.~\ref{fig:area_map_eps}).

\begin{figure*}
\includegraphics[width=\textwidth]{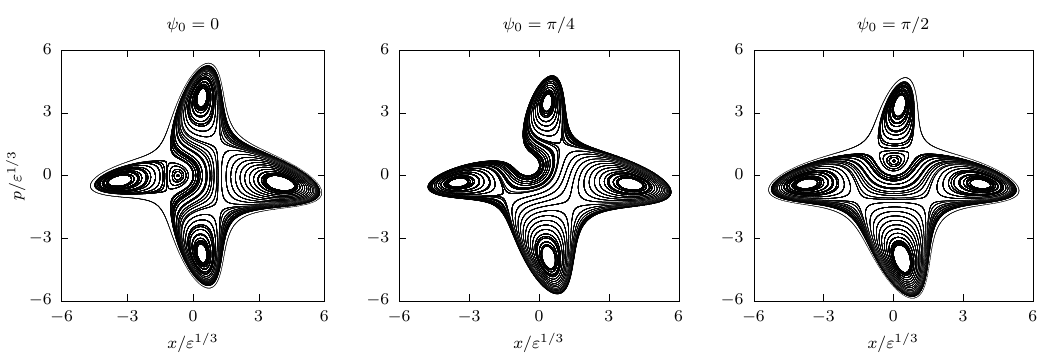}
\caption{Phase space portraits of the Poincaré map of Eq.~\eqref{eq:doubleresmap} sampled every $4$ iterations for three values of the initial phase $\psi_0$. The phase of exciter selects the island that becomes larger than the others (the parameters values are $\omegar/(2 \pi)=1/4$, $\Delta=0$, $\eps=10^{-4}$, $\kappa=\num{0.1}$, $\delta/\eps^{2/3}=0.75$).}
\label{fig:phsp_map_ph}
\end{figure*}

The phase-space topology is also influenced by the exciter phase $\psi_0$. It is straightforward to verify that if $\psi_0=m \, \pi/2, \,\, m \in \mathbb{N}$, the resulting map coincides with the case $\psi_0=0$ after $m$ additional turns. Therefore, the effect of the exciter phase is a rotation of the phase space. When $\psi_0$ assumes intermediate values between $m\pi/2$ and $(m+1)\pi/2$, the result is more interesting. In fact, as $\psi_0$ increases, the area of the main island decreases and that of an adjacent island increases. At the middle point, when $\psi_0= (2m+1)\pi/4$, two islands of equal size become dominant. This transition is clearly visible in Fig.~\ref{fig:phsp_map_ph}.

\subsection{Hamiltonian model}
\label{sec:ham}

We introduce a Hamiltonian model to study the phase-space structure of the map model~\eqref{eq:doubleresmap} when $\Delta_n=0$ and $\omegar/(2 \pi)=1/4$. We take advantage from the Birkhoff Normal Form expansion of the Hénon map close to the fourth-order resonance up to the third order, which provides an interpolating Hamiltonian of the form 
\begin{equation} 
\begin{split}
\ham(\phi,J) & = \omega_0(\lambda) J + J^2\qty[\frac{\Omega_2}{2} + A\cos(4\phi-2\pi t)] \, , 
\label{eq:mainH0}
\end{split}
\end{equation}
where $\omega_0(\lambda) = \pi/2+\delta(\lambda)$, $\lambda=\epsilon t$, and the Normal Form computation gives the coefficients of the Hamiltonian~\cite{yellowreport}, \ie
\begin{equation}
\begin{aligned}
\Omega_2 &= -\frac{1}{8}-\frac{3}{8}\kappa + \mathcal{O}(\delta),\\
A &= \frac{1}{16}(1-\kappa)+\mathcal{O}(\delta)\,, 
\end{aligned}
\end{equation}

Then, we introduce the resonant exciter that corresponds to a linear term $\propto x$ in the Hamiltonian according to the definition in the map~\eqref{eq:doubleresmap}. To express $x$ as a function of action $J$ and angle $\phi$, we introduce the complex variables $z=x-ip$, its complex conjugate $z^*$, the conjugating function of the Normal Form transformation $\bm\Phi(\zeta,\zeta^*)$, and action $2J=\zeta\zeta^*$. At the first perturbation order in $J$ the expansion of $x$ reads
\begin{equation}
x = \frac{z+z^*}{2} = \frac{\zeta + \zeta^*}{2} + \frac{\Phi(\zeta,\zeta^*) + \Phi^*(\zeta,\zeta^*)}{2}\,,
\label{eq:dres_q}
\end{equation}
with 
\begin{equation}\Phi(\zeta,\zeta^*) = \sum_{n\ge 2}\sum_{k=0}^n \phi_{k,n-k}\zeta^k\zeta^{*n-k}\,.
\end{equation}

In action-angle variables, Eq.~\eqref{eq:dres_q} becomes 
\begin{equation}
\begin{split}
    x(\phi,J) &= \sqrt{2J}\cos\phi \\ & \qquad + \sum_{n\ge 2}(2J)^{n/2}\sum_{k\le n}\Re(\phi_{k,n-k})\cos((2k-n)\phi)\\ &=\sum_{\ell\ge 1} x_\ell(J)\cos(\ell \phi) \,,
\end{split}
\end{equation}
where
\begin{equation}
\begin{split}
    x_1(J) &= \sqrt{2J}\qty[ 1 + \sum_{n\ge 1} (2J)^n\Re(\phi_{n+1,n}+\phi_{n,n+1})],\\
    & \\
    x_\ell(J) &= \sum_{n}(2J)^{n/2} \Re\qty(\phi_{\frac{n+\ell}{2},\frac{n-\ell}{2}}+\phi_{\frac{n-\ell}{2},\frac{n+\ell}{2}})\qquad \ell > 1\,,
\end{split}
\end{equation}
and it is worth noting that the numbers $\ell$ and $n$ in the previous equation must have the same parity.
%
%
%
%

The Hamiltonian of Eq.~\eqref{eq:mainH0}, with the contribution of the exciter, can be then written as
\begin{equation}
\begin{split}
\ham = \omega_0(\lambda) J &+ J^2\qty[\frac{\Omega_2}{2} + A\cos(4\phi-2\pi t)] + \\
    & + \eps \sum_\ell x_\ell(J)\cos(\ell\phi)\cos(\omega t + \psi_0)\,.
\end{split}
\label{eq:ham_dres1}
\end{equation}
The last term of Eq.~\eqref{eq:ham_dres1}, can be rewritten introducing the slow angle $\theta=\phi - \pi t/2$ according to
\begin{equation}
\begin{split}
\cos(\ell\phi)&\cos(\omega t+\psi_0) = \\  &= \frac{1}{2}\left[\cos(\ell\theta + \frac{\ell+1}{2}\pi t + \Delta(\lambda) t + \psi_0) + \right. \\ & \qquad \left. \cos(\ell\theta + \frac{\ell-1}{2}\pi t - \Delta(\lambda) t - \psi_0)\right] \, .
\end{split}
\label{eq:dres_avg}
\end{equation}
where $\Delta (\lambda)=\omega-\omega_r$. When $\Delta(\lambda)=0$, the time average of Eq.~\eqref{eq:dres_avg} is $\cos(\theta+\psi_0)/2$ only for $\ell=1$, and zero otherwise, and the averaged expansion of $x$, up to order $J^3$, becomes
\begin{equation}
\begin{split}
    \av{x(J)} &= \sqrt{2J}\qty[1 + 2J\,\Re(\phi_{21}+\phi_{12})]\\
              &= \sqrt{2J}\qty(1+c_1 J)\,,
\end{split}
\end{equation}
where $c_1$ is a constant term whose value can be retrieved from the computation of the terms of $\Phi$ in normal form (in our case, we have $c_1=2 \Re(\phi_{21}+\phi_{12}) = \mathcal{O}(\delta)$). Neglecting terms of order $\mathcal{O}(\delta)$ in the resonant term, we have the Hamiltonian
\begin{equation}
\begin{split}
\ham = \delta(\lambda) J & + J^2\qty(\frac{\Omega_2}{2} + A\cos 4\theta) + \\
& +\frac{1}{2}\eps\sqrt{2J}\cos(\theta + \psi_0)\, .
\end{split}
\label{eq:ham_dres_fin}
\end{equation}

\begin{figure*}
\includegraphics[width=\textwidth]{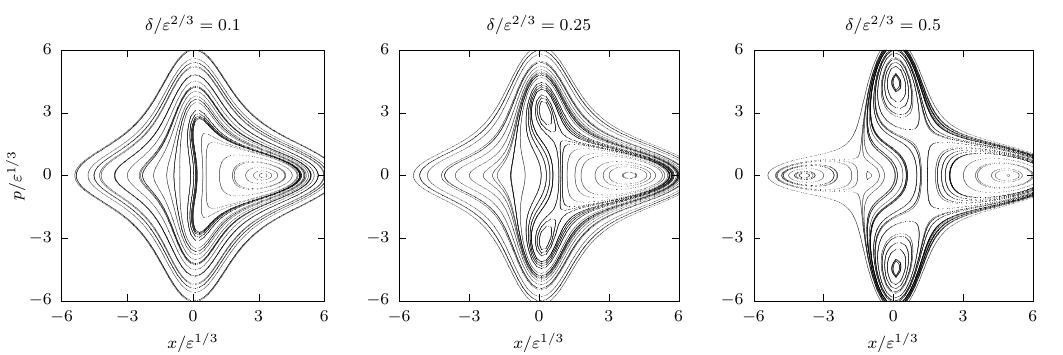}
\includegraphics[width=\textwidth]{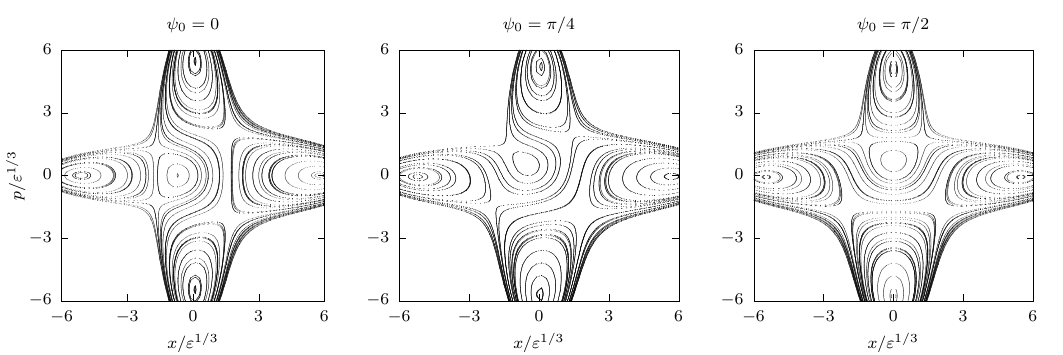}
\caption{Top: Phase-space portraits of the Hamiltonian of Eq.~\eqref{eq:ham_dres_fin} with $\psi_0=0$ and $\kappa=0.1$, for three values of the scaled parameter $\hat\delta=\delta/\eps^{2/3}$ representing three topological structures, which correspond to those shown in Fig.~\ref{fig:phasespace_map} for the map model. The coordinates have also been re-scaled by the factor $\eps^{1/3}$. Note that for $\kappa=0.1$, the two critical values of $\hat \delta$ where the solutions bifurcate are $\hat\delta_1 = \num{0.348}$ and $\hat\delta_2 = \num{0.179}.$ Bottom: Phase-space portraits of the Hamiltonian of Eq.~\eqref{eq:ham_dres_fin} with $\hat \delta=0.75$ and $\kappa=0.1$, for three values of the initial phase $\psi_0$, corresponding to Fig.~\ref{fig:phsp_map_ph} for the map model. The same coordinate scaling has been used.}
\label{fig:phsp_ham}
\end{figure*}

In Fig.~\ref{fig:phsp_ham}, we observe that, depending on the values of $\delta$, $\eps$ and $\psi_0$, the phase-space portraits of Eq.~\eqref{eq:ham_dres_fin} present the same features as those of Eq.~\eqref{eq:doubleresmap}, with the appearance of resonance islands of unequal area. However, it is worth noting that the position of the fixed points in the phase space is not the same for the map or for the Hamiltonian. This is due to the difference in the dependence of the orbit frequency on the action that, in the Hamiltonian, is truncated at the second order in $J$. The presence of a positive cubic term, \ie $\kappa >0$, in the map is necessary to have closed separatrices in the phase space of the Normal Form Hamiltonian at order $O(J^2)$. Resonance islands can also be found for $\kappa \le 0$, but to retrieve the same topology, one must perform an additional step in the Normal Form computation, which makes the model overly complicated.

From Eq.~\eqref{eq:ham_dres_fin} we can also retrieve some scaling properties of the parameters. Let $J = \hat J \eps^{2/3}$ and $\delta=\hat \delta \eps^{2/3}$, the Hamiltonian can be scaled by defining a slow time $\hat t= \eps^{4/3} t$, and we obtain
\begin{equation}
\begin{split}
\hat\ham = \hat\delta (\lambda) \hat J & + \hat J^2\qty( \frac{\Omega}{2} + A\cos 4\theta) + \\
& + \frac{1}{2}\sqrt{2\hat J}\cos(\theta-\psi_0) \, , 
\end{split}
\label{eq:ham_dres_fin_scale}
\end{equation}
which means that if the ratio $\delta/\eps^{2/3}$ is kept constant, the resulting phase-space portraits of the Hamiltonian of Eq.~\eqref{eq:ham_dres_fin} are similar (up to a scaling of the action by a factor $\eps^{-2/3}$). From this consideration one observes that, since $\eps$ and $\delta$ are found in reciprocal positions, the phase-space topologies that the system crosses, increasing $\eps$ at constant $\delta$, have a reversed order when increasing $\delta$ at constant $\eps$.

Appendix~\ref{app:fps} discusses in detail the fixed points and resonance islands of the Hamiltonian~\eqref{eq:ham_dres_fin}. In particular, we establish the dependence of the phase-space topology on the parameter $\hat\delta = \delta/\eps^{2/3}$, and we identify two critical values of $\hat\delta$ associated with the emergence of fixed points: $\hat\delta_1 = 3\kappa^{1/3}/4$ and $\hat\delta_2=(3/4)\qty[\kappa^{2/3}(\kappa+1)^{1/3}(\kappa-1)^{-2/3}]$.

The relevance of these concepts can be seen in Fig.~\ref{fig:phsp_ham}, where one observes the nice correspondence between the Hamiltonian Normal Form and the map phase-space portraits of Figs.~\ref{fig:phasespace_map} and ~\ref{fig:phsp_map_ph}, when the same values of $\delta/\eps^{2/3}$ and $\psi_0$ are used. The Hamiltonian description succeeds in reproducing all possible phase-space topologies of the map depending on the scaled parameter $\hat\delta$ and on the initial exciter phase.

\section{Qualitative analysis of the trapping process} \label{sec:trapping}

Consider the evolution of an orbit of the map defined by Eq.~\eqref{eq:doubleresmap}. At the beginning, \ie $n=0$, we set $\omega_{0,0}=\pi/2 - \delta$ and $\eps_0=0$, and then the map is iterated for $2N$ turns. During the first $N$ turns, $\eps$ increases linearly from $0$ to $\eps_N=\eps$, while $\omega_0$ remains constant. During the second $N$ turns, we keep $\eps_n=\eps$ and linearly increase $\omega_0$ by $2\delta$. With this strategy, the first part of the process is meant to slowly introduce the effect of the exciter, whereas in the second part the actual resonance-crossing process takes place. 

\begin{figure}
\centering
\includegraphics[width=0.5\columnwidth]{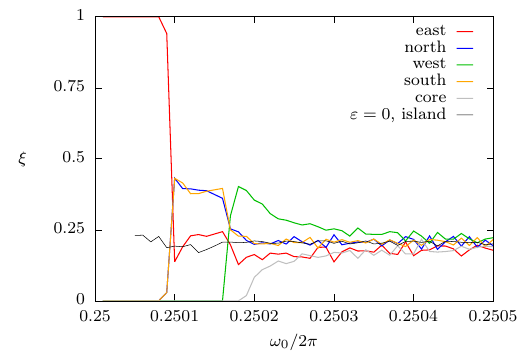}
\caption{Theoretical probability $\xi$, according to Eq.~\eqref{eq:neish}, for a particle to be trapped in an island or the core depending on the value of the map main frequency $\omega_0$. For the sake of comparison with the case of the Hénon map without exciter, the continuous black line represents the trapping probability in an island when $\eps=0$. This plot has been made by numerically computing the derivative of the area of each phase-space region of the map of Eq.~\eqref{eq:doubleresmap} (parameters values: $\Delta=0$, $\eps=10^{-4}$, $\psi_0=0$).}
    \label{fig:areadiff}
\end{figure}

To understand the resonance trapping process in terms of the separatrix-crossing theory developed by Neishtadt, whose main result is the probability formula given in Eq.~\eqref{eq:neish}, we consider the dependence of the values of the area of the islands and the core region on the parameter $\hat \delta$ (see Fig.~\ref{fig:area_map_eps}), and the dependence of the trapping probability $\xi$ of an orbit in any of the phase-space region as a function of $\omega_0$, \ie of $\delta$, which varies linearly (see Fig.~\ref{fig:areadiff}), which are related to the process we are going to describe in this section.

Depending on the values of $\eps$ and $\delta$, during the first phase of the process, resonance islands may or may not be present in the phase space. However, as the frequency $\omega_{0, n}$ increases and goes beyond the resonant value $\omegar$, the four islands eventually appear. 

If we follow the evolution of the phase-space topology of the $4$th iterate of the map~\eqref{eq:doubleresmap} during the trapping process, it is possible to understand how a given initial condition will be trapped in one of the possible phase-space regions. In the initial state of the modulation process, when $\delta<0$ and $\eps=0$, there are no resonance islands, the stable fixed point is at the origin of the phase space, and inside the stability domain of the map, the orbits of the initial conditions rotate around the stable fixed point.

As $\eps$ increases, the central fixed point moves in the $x$ positive direction (or, using our `geographical terminology, it moves eastward). As this process is adiabatic, the orbits moving around the central fixed point will also be displaced eastward, and when $\eps$ reaches the final value, the orbits will rotate centered around the new position of the stable fixed point. All initial conditions at low amplitude will remain in the stable basin of the fixed point, which defines the main island. 

For $\hat\delta > \hat\delta_2$ two new islands appear in the phase space: One is the North island and the second is the South island. As can be seen from the evolution of the areas $A_i$ (see Fig.~\ref{fig:area_map_eps}), the size of the secondary islands increases at the same rate as that of the principal, but their area is always smaller. Therefore, as Fig.~\ref{fig:areadiff} shows, there is a range of $\omega_0$ (\ie of $\delta$), for which the area derivative \wrt $\delta$ of these secondary islands is higher than that of the principal one. Hence, in this range the particles will preferentially be trapped in the North or South islands. As the area of the island structure grows while $\omega_0$ is changed, orbits whose initial conditions are in the outer region will be initially trapped in the East island and afterward also in the new islands, and we have a minimum value for the action to be trapped in them.

Finally, when $\hat\delta > \hat\delta_1$, a new separatrix appears and the last island is created in the phase space on the opposite side of the principal one. This is the smallest island, but also has a peak of the time derivative of its area, which corresponds to a peak in the probability of trapping. Therefore, a certain fraction of high-amplitude particles will enter the last one. Together with the smaller island, a new fixed point is created in the central region, and the core area will also trap some high-amplitude particles. Finally, for even higher values of $\delta$ we see that the area derivatives of the five regions tend to be equal and, therefore, no island is privileged when $\hat\delta$ is sufficiently large. Note that this also explains the absence of dominant islands for the standard Hénon map, corresponding to $\eps=0$.

Figure~\ref{fig:areadiff} also reports the trapping probability for the case of the nominal Hénon map, \ie without an external exciter. A single curve is shown as the four islands behave the same. Furthermore, note that according to Ref.~\cite{yellowreport}, for $\eps=0$ one should expect $\xi \to 1/4$ if $\omega_0/2\pi \to 0.25$, but the numerical computation of $\xi$ is made difficult by the small values of the islands and core areas close to the resonance, which explains the missing part of the curve.

\section{Simulation results} \label{sec:simresults}

Detailed numerical simulations have been carried out with the aim of understanding and controlling the trapping process and adiabatic transport inside the stable islands for the map~\eqref{eq:doubleresmap}. The quantities $\eps$ and $\delta$ have been linearly varied as a function of the number of turns to observe how the trapping in each island depends on various parameters of the system.

We study the evolution of $N_\mathrm{p}=\num{3e3}$ particles, distributed with an initial action average $\av{J_0} = \av{x_0^2+p_0^2}/2$, after $2N$ turns of the map, using the process described in Section~\ref{sec:trapping}. At the end of the modulation, each particle orbit can be classified using the main and secondary tune values. Tune analysis is performed to provide an accurate identification of the region in which each initial condition is trapped at the end of the evolution process. The region identification is performed by looking at the frequency of the final orbits, according to the approach described in the Appendix~\ref{app:tune}.
\begin{figure*}[htb]   
\includegraphics[width=\textwidth]{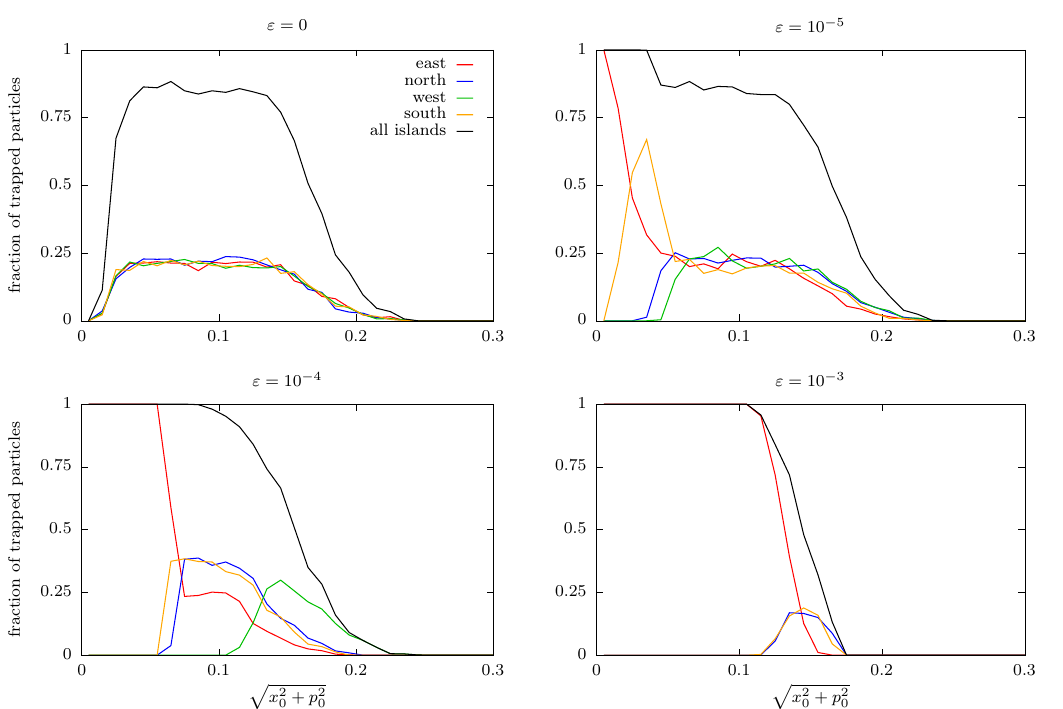}
    \caption{Fraction of particles trapped in each island and their sum for initial conditions of the map of Eq.~\eqref{eq:doubleresmap} as a function of the initial amplitude $\sqrt{2J_0}=\sqrt{x_0^2 + p_0^2}$ of the particles for four values of exciter strength $\eps$. Data for the graphs were collected by simulating several uniform initial annular distributions, each containing $N_\mathrm{p}=\num{1e3}$ particles, for each bin of the histogram (bin width $0.01$). Note that the upper-left plot corresponds to a standard Hénon map without external exciter (parameters values:  $\omegar/(2\pi)=1/4$, $\delta/(2\pi)=\num{5e-4}$, $\kappa=0.1$, $\psi_0=0$, $N=10^5$).}
    \label{fig:plot_distr_map}
\end{figure*}
\begin{figure*}\includegraphics[width=\textwidth]{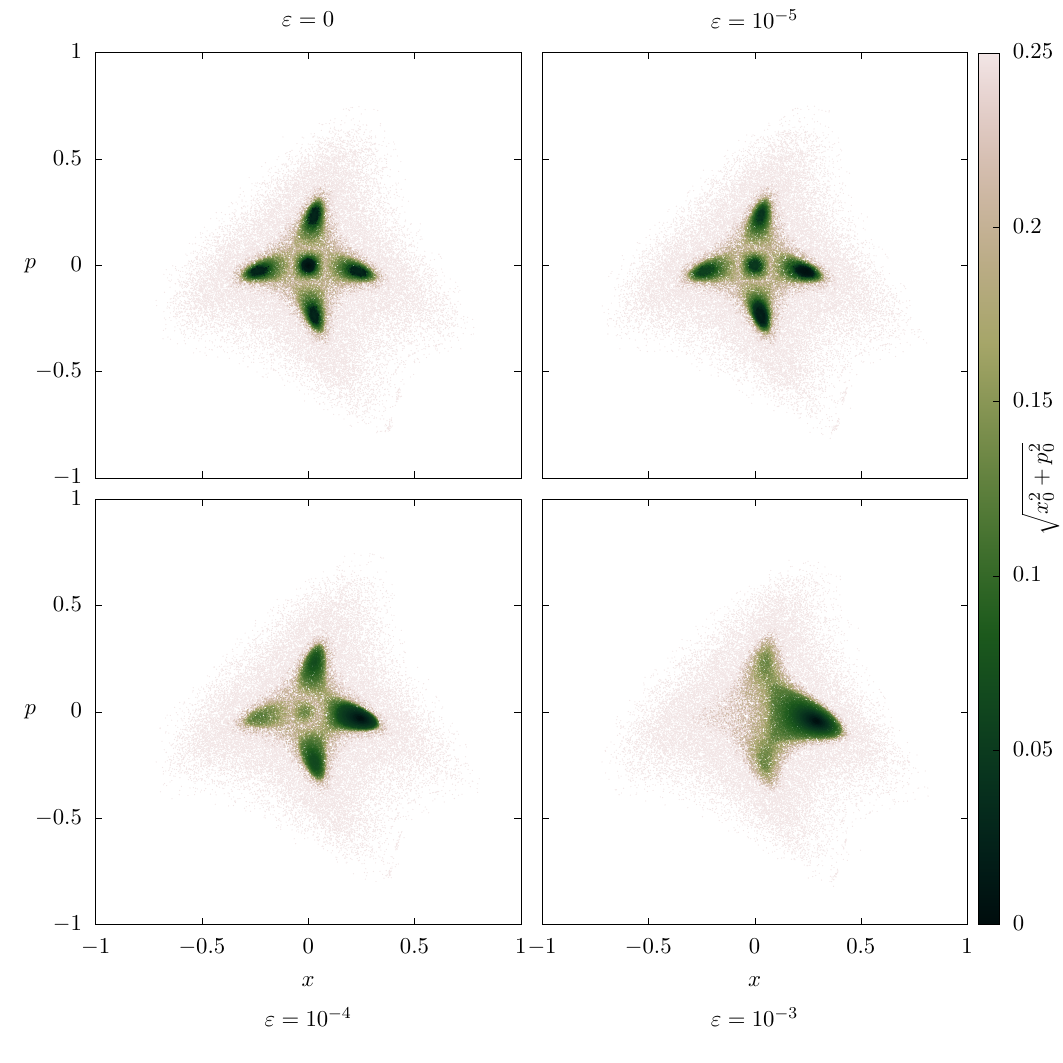}
    \caption{Final state in the phase space $(x,p)$ for a family of initial annular uniform distributions of particles with amplitudes $0\le \sqrt{x_0^2 + p_0^2} \le 0.25$ for four values of $\eps$ using the same parameters of the numerical simulations whose results are shown in Fig.~\ref{fig:plot_distr_map}. Note that the upper-left plot corresponds to a standard Hénon map without external exciter. The color scale encodes the initial amplitude.}
    \label{fig:plot_distrxy_map}
\end{figure*}

The first study consists of scrutinizing how the initial distribution of particles is transformed into the final one, \ie how the initial conditions are shared between the islands. Figure~\ref{fig:plot_distr_map} shows the fraction of particles found in each island, for a family of uniform annular initial distributions, as a function of their amplitude $\sqrt{x_0^2 + p_0^2}=\sqrt{2J_0}$. For each histogram bin, $N_\mathrm{p}=\num{1e3}$ particles have been generated, uniformly distributed in an action interval of width of $0.01$, and we calculated, at the end of the modulation process, how many particles of each annulus were trapped in each island. This is repeated for four values of the final exciter amplitude. We observe that when $\eps=0$, all islands have the same behavior and capture particles at the same amplitudes. Furthermore, at small amplitude, the islands do not trap any initial condition, and these are then left orbiting around the stable fixed point at the origin of the phase space. As $\eps$ increases, we see that particles, even at low amplitudes, are trapped more and more in the main island (the East island), whereas smaller islands begin to trap particles only at higher amplitudes. This behavior is not possible for the standard Hénon map, which is shown in the top-left plot, and is the main effect that led to the use of an exciter in the application of this process to accelerator physics. 

\begin{figure*}
\includegraphics[width=0.95\textwidth]{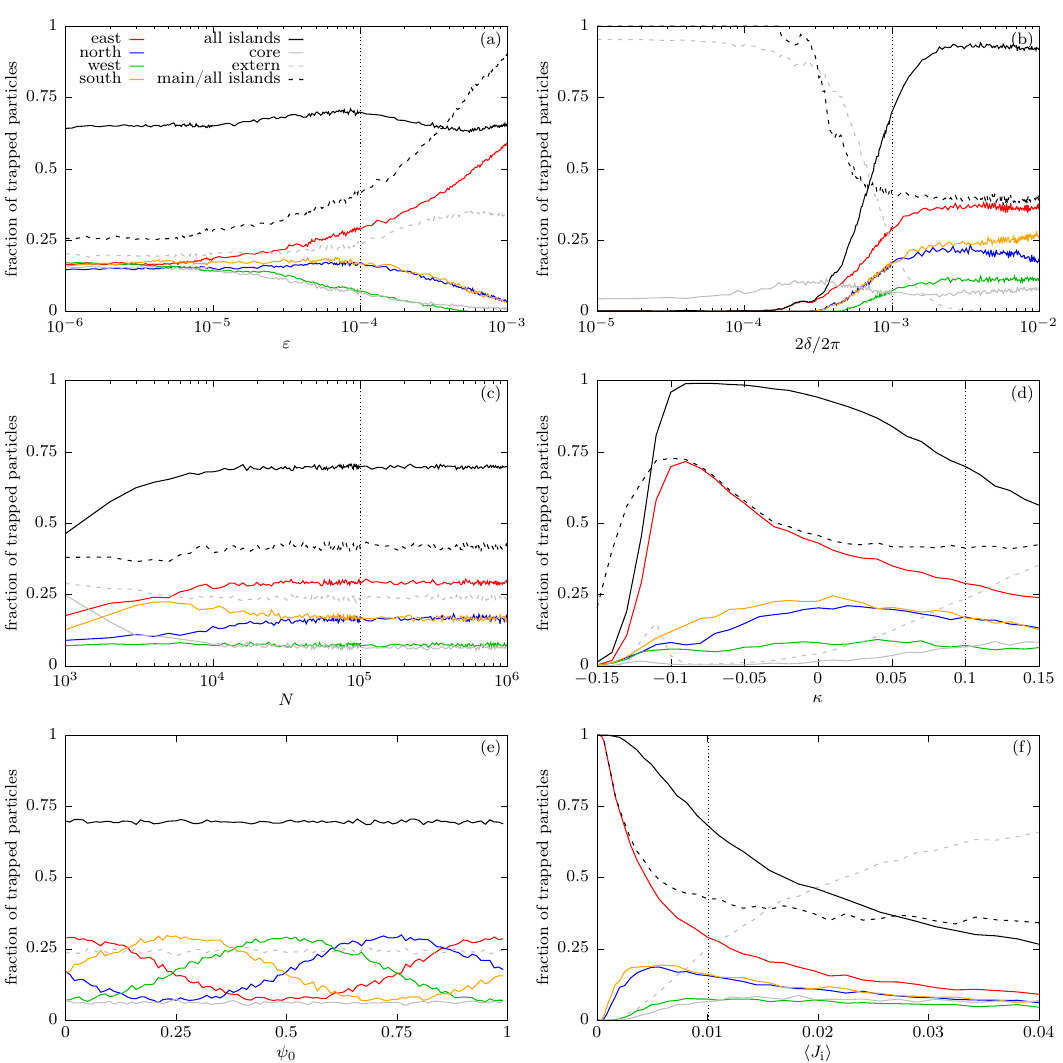}
    \caption{Final trapping fraction in each island, in the core region, in the external region, and ratio between the final trapping fraction in the main island and that of the sum of all islands for a normal random distribution of $N_\mathrm{p}=\num{3e3}$ initial conditions under the evolution of the map of Eq.~\eqref{eq:doubleresmap} with $\omega/(2 \pi)=1/4$ changing one parameter at a time and keeping all the others fixed to the default values of $\eps=10^{-4}$, $\delta/(2 \pi)=\num{5e-4}$, $\kappa=0.1$, $\psi_0=0$, $\av{J_\text{i}}=0.01$, $N=10^5$ (the vertical line represents the default value of the parameter). Starting from the top-left graph, the final trapping fraction of the distribution is plotted against the exciter strength value $\eps$ (log scale), the main frequency change $2\delta$ (log scale), the number of time steps $N$ (log scale), the strength of the cubic nonlinearity $\kappa$ (linear scale), the initial phase of the exciter $\psi_0$ (linear scale) and the average initial action $\av{J_\text{i}}$ (linear scale). The legend shown on the top-left graph is valid for all other graphs.}
    \label{fig:plot_phase_map}
    \label{fig:plot_eps_map}
    \label{fig:plot_delta_map}
    \label{fig:plot_oct_map}
    \label{fig:plot_emitt_map}
    \label{fig:plot_n_map}
\end{figure*}

The final distributions of the particles are shown in Fig.~\ref{fig:plot_distrxy_map}, where the color scale encodes the initial amplitude. As the value of $\eps$ increases, the asymmetry between the islands becomes more and more visible. Furthermore, the set of particles left at the center of the phase space is reduced until it is completely emptied at the largest value considered, \ie $\eps=10^{-3}$. It should be noted that low-amplitude particles are trapped in the inner region of the islands, or in the center, for the case $\eps=0$. However, as $\eps$ increases, the low-amplitude particles are trapped in the main island, while the other islands trap the higher-amplitude particles, in agreement with the results shown in Fig.~\ref{fig:plot_distrxy_map}. The red halo visible in Fig.~\ref{fig:plot_distrxy_map} is generated by large-amplitude particles that are not trapped in any island. 

Several additional studies were carried out with the aim of assessing the dependence of the trapping process on the parameters of the model considered. We defined a set of parameter values that we consider as the default, when others are varied to study the dependence of the trapping fraction. The default values consist of the maximum exciter amplitude ($\eps=\num{1e-4}$), its constant frequency ($\omega/(2\pi)=1/4$), its phase ($\psi_0=0$), and we also include the octupolar coefficient ($\kappa=0.1$). The rotation frequency ranges from $\omega_{0,0}/(2\pi)=\omega_{0,N}/(2\pi)=1/4-\delta/(2\pi)$ to $\omega_{0,2N}/(2\pi)=1/4+\delta/(2\pi)$, with $\delta/(2\pi)=\num{5e-4}$. The total number of simulated turns is $2N=\num{2e5}$, and the initial conditions are normally distributed, inspired by what occurs in accelerator physics applications, in both $x$ and $p$, with standard deviations $\sigma_x=\sigma_p=0.1$, which account for an initial average action of $\av{J_0}=0.01$. Although the approach followed here provides an essential understanding of the main features of the process under study as a function of the various model parameters, it goes without saying that a global optimization of the entire process \eg in the sense of defining the target trapping fraction in each island, which would be needed in the case of realistic applications, would require a global multi-parameter optimization, not a scan of a single parameter at a time.

Figure~\ref{fig:plot_eps_map}~(a) shows the fraction of initial conditions trapped in each island and the ratio between the trapping in the main island (as $\psi_0=0$, this is the East island) as a function of the final exciter amplitude $\eps$. We observe that the difference in trapping between islands increases with $\eps$. The peculiar behavior of the dominant island is clearly visible: For all values of $\eps >10^{-5}$, it can be observed that the East island captures the highest fraction of initial conditions (up to $\approx 60\%$ of the initial distribution, which represents $\approx 90\%$ of the islands), followed by the North and South islands and then the smallest West island. The effect is more prominent for large values of $\eps$, whereas for $\eps\to 0$ all islands tend to behave the same. Note that the situation at $\eps=10^{-6}$ is almost indistinguishable from that of the Hénon map without exciter. However, the total number of particles in the islands increases only slightly, and until $\eps\sim \num{1e-4}$. Above that value, the total number of particles in the islands is reduced because a larger number of particles are expelled to the external region of the phase space.

If we study the trapping properties as a function of the frequency excursion $\delta$ (for which $\omega$ varies from $\omega-\delta$ to $\omega+\delta$), (see Fig.~\ref{fig:plot_eps_map}~(b)), we observe that if the frequency variation is not large enough, most particles will not be trapped in the islands, but will be expelled in the external region. In this situation, the area of the island structure is not large enough to match the initial action of most particles, which therefore do not get trapped. At the other extreme, for large values of $2 \delta/(2\pi)$, a small decrease in trapping is observed in some islands, which is due to the loss of particles due to escape to infinity induced by nonlinearities. It is worth observing and stressing how the parameter $\delta$ can be used to improve the trapping fraction with respect to the default case.

The role of adiabaticity is visible in the results shown in Fig.~\ref{fig:plot_n_map}~(c), where we present the trapping fraction in each region as a function of the number of turns $N$ during which the modulations are performed. The fraction of particles successfully trapped increases with $N$ reaching a plateau at $N\sim 10^4$, therefore confirming that the default value of $10^5$ represents a sufficiently slow modulation to apply adiabatic theory to the system.

The role of the cubic nonlinearity present in the map is studied and visible in the results shown in Fig.~\ref{fig:plot_oct_map}~(d), where we analyze the trapping fractions as a function of $\kappa$. The main effect of the parameter $\kappa$, which represents the strength of the cubic nonlinearity, is to cause a deformation of the phase space, changing the shape of the islands. Furthermore, it also reduces the extent of the stable region of the phase space, \ie the region where a bounded motion occurs, causing escape to infinity of orbits of several initial conditions. This behavior becomes dominant in the region $\kappa < -0.1$. However, for large positive values of $\kappa$, the islands become smaller and the time derivative of their surface changes: Fewer particles are trapped in the islands, although the ratio between the trapped fraction in the main island and all islands does not change much from what was observed when $\kappa=0$.

Figure~\ref{fig:plot_phase_map}~(e) shows the fraction of initial conditions trapped in each island as a function of the exciter phase $\psi_0$. We see that, depending on $\psi_0$, a different island will capture most of the initial conditions, as we observed when looking at the phase space of the stroboscopic map with different exciter phases. This change is a smooth function of $\psi_0$, and there are special values of $\psi_0$ for which two dominant islands are present and whose size is comparable. The trapped fraction in all islands, as well as the final population of the core and of the external regions, is constant \wrt $\psi_0$.

Finally, Fig.~\ref{fig:plot_emitt_map}~(f) shows the evolution of the trapping fraction in each phase-space region as a function of the action average of the initial distribution $\av{J_\text{i}}=(\sigma_x^2 + \sigma_p^2)/2$. This plot represents an integration over the initial normal random distribution of what was observed in Fig.~\ref{fig:plot_distr_map}, which was obtained using concentric annular distributions of the initial conditions. As particles close to the origin are more likely to be trapped in the main island, this effect is more prominent in the case of initial distributions corresponding to small values of $\av{J_\text{i}}$. However, the other islands are populated with initial conditions starting at higher amplitudes. Therefore, wider distributions present a more balanced output, reducing the ratio between the trapped fraction in the main island and the total fraction trapped in all the islands. This improved balance is achieved at the price of a reduced total number of trapped particles in the islands, as the high-amplitude fraction of the distribution of initial conditions lies beyond the phase-space region covered by the islands. In this case, also, a higher value of the parameter $\delta$ would allow for an optimization of the trapping process.

\begin{figure}[htb]
    \centering
    \includegraphics{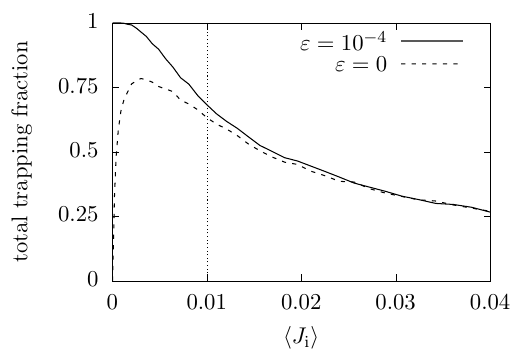}
    \caption{Comparison between the total trapping fraction into the four islands using a double-resonant approach ($\eps=10^{-4}$, black line), and the Hénon map (with cubic nonlinearity, but no external exciter) splitting ($\eps=0$, dashed line) as a function of the action average of the normal distribution of the initial conditions. The vertical dotted line represents the value of $\av{J_\text{i}}$ used in other studies (parameters values: $\omegar/(2 \pi)=1/4$, $\delta/(2 \pi)=\num{5e-4}$, $\kappa=0.1$, $\psi_0=0$, $N=10^5$, $N_\mathrm{p}=\num{3e3}$).}
    \label{fig:compare_henon}
\end{figure}

As a last consideration, we would like to highlight the peculiar behavior introduced by the double-resonance condition, which is shown in Fig.~\ref{fig:compare_henon}, as this can play a fundamental role for applications. In that plot, the fraction of initial conditions trapped in the islands is shown as a function of the average action of the normal distribution of the initial conditions. Two cases are reported, namely the Hénon map (with cubic nonlinearity but no external exciter) and the special map with external exciter fulfilling a double-resonance condition.

The difference in behavior between the two models is clearly seen for initial distributions with an average action value close to zero. In this condition, the Hénon map features a total fraction of initial conditions trapped in the islands that goes to zero, while when the exciter is in action, the fraction of initial conditions trapped in the islands tends to one. The difference between the two models fades away as $\av{J_\mathrm{i}}$ increases. Therefore, the use of the exciter with a double-resonance condition allows efficient control of particle trapping even for the case of initial distributions of a rather small extent in phase space, for which the natural trapping would naturally be very low.
\section{Conclusions} \label{sec:conc}
In this paper, we have presented the detailed study of a double-resonance condition for a 2D model of a nonlinear Hénon-like map. The double resonance occurs because the main frequency of the map is close to a rational value ($1:4$ in our case), and the corresponding resonance is excited by the nonlinearities of the map and the frequency of a time-dependent external exciter, whose frequency is in resonance $1:1$ with the main frequency of the system. 

The corresponding Normal-Form Hamiltonian model has been constructed and studied in detail. This model is of paramount importance for understanding the scaling properties for the various model parameters and the fixed point and islands topology in the phase space, which can be transferred to the map model. The Hamiltonian model clearly indicated that the key feature introduced by the double-resonance condition is an asymmetry in the island structure of the Henon map. In fact, one island grows and becomes dominant with respect to the other three. Furthermore, the structure of the single separatrix of the Hénon map is completely altered and split into separatrices. These effects provide an original and new phase-space foliation that has deep implications for adiabatic trapping and transport in phase space.

After the analysis of the phase-space topology, the study of adiabatic trapping in the islands was carried out. According to the theory of adiabatic resonance crossing, a distribution of initial conditions can be partially trapped into the various structures in phase space, generated by separatrices of the frozen system. By means of detailed numerical simulations, we showed that trapping actually occurs and studied the dependence of its properties on the various model parameters. These studies can be used to optimize the trapping process, which is an essential aspect for an application of this novel double-resonance system. It is important to stress that the obtained results agree with the experimental observations and that the use of an external exciter satisfying a double-resonance condition allows the trapping probability to increase in the islands. 

\section*{Acknowledgments}
We are indebted to Prof. A.~Neishtadt for several discussions and useful suggestions. 

\appendix 

\section{Fixed points and resonance islands} \label{app:fps}

To describe the properties of the Hamiltonian phase space (see Eq.~\eqref{eq:ham_dres_fin}) and to discuss the resonance trapping phenomena when the parameters $\eps$ and $\omega_0$ are varied, we start from the fixed-point analysis in case $\psi_0=0$. Since
\begin{equation}
    \pdv{\ham}{\theta} = -4AJ^2\sin4\theta  -\frac{\eps}{2}\sqrt{2J}\sin\theta \, , 
\end{equation}
we have exact solutions to the fixed-point equation for $\theta=0,\pi$ and approximate solutions at $\mathcal{O}(\eps)$, close to $\theta=\pm \pi/4$ and $\theta=\pm \pi/2$. In the first case, it is more convenient to write the Hamiltonian using Cartesian coordinates $X=\sqrt{2J}\cos\theta,\,Y=\sqrt{2J}\sin\theta$, and consider the solutions of $\pdv*{\ham}{X}$ at $Y=0$ that reduces to the cubic equation
\begin{equation}
    -\frac{\kappa}{3}X^3 + \delta X + \frac{\eps}{2} = 0 \, .
    \label{eq:cubic}
\end{equation}
The exact solutions could be retrieved using the well-known Cardano formulae. However, for a qualitative analysis of the existence of solutions, it suffices to study the sign of the discriminant $\kappa\delta^3 - 27\kappa^2\eps^2/64$. Assuming $\kappa>0$ and $\eps>0$, we have three real solutions if 
\begin{equation}
 \frac{\delta}{\eps^{2/3}} > \frac{3}{4}\kappa^{1/3} = \hat\delta_1   
\end{equation}
and a single real (positive) solutions otherwise. According to Descartes' rule of signs, one of the three real solutions is found on the positive $X$ semiaxis (i.e., at $\theta=0$), and two on the negative one ($\theta=\pi$). 
The three real solutions are the elliptic points of three stable regions: The center and two islands, the West and East islands, according to our nomenclature. For $\delta$ smaller than the critical value, the real single solution coincides with the origin of the phase space when $\eps=0$, and represents the displacement of the central region, which is found at a distance $\approx \eps$ for small values of $\eps$ and $\approx \eps^{1/3}$ for large values of the exciter amplitude.
The solutions for $Y=0$ are given by the equation
\begin{equation}
    -\frac{\kappa^2-1}{4\kappa} X^3 + \delta\frac{\kappa -1}{\kappa}X + \frac{\eps}{2} = 0
    \label{eq:cubic2}
\end{equation}
with
\begin{equation}
    Y = \pm\sqrt{\frac{4\delta}{\kappa} - \frac{X^2}{\kappa}} \, , 
\end{equation}
which gives up to six extra fixed points, symmetrical \wrt the axis $X=0$.

Assuming $\kappa < 1$, Eq.~\eqref{eq:cubic2} has up to three real roots if
\begin{equation}
    \frac{\delta}{\eps^{2/3}} > \frac{3}{4}\qty[\frac{\kappa^2(\kappa+1)}{(\kappa-1)^2}]^{1/3} = \hat\delta_2
\end{equation}
and one real root otherwise. For $\kappa<1/3$ we have $\hat\delta_2<\hat\delta_1$.

In the top three plots of Fig.~\ref{fig:phsp_ham} we show three different phase-space portraits of the Hamiltonian of Eq.~\eqref{eq:ham_dres_fin} for $\psi=0$, $\kappa=0.1$ and three possible configurations depending on the value of $\hat\delta=\delta/\eps^{2/3}$, (where the scales of $x$ and $p$ depend on the value of $\eps$). For $\hat\delta<\hat\delta_2$, we only have a fixed point on the positive semiaxis $x$, which means that a single island (the East island) is present in the phase space. For $\hat\delta_2<\hat\delta<\hat\delta_1$ we still have a fixed point on the positive $x$ semiaxis and the other four symmetric fixed points (two elliptic and two hyperbolic), which create two equal islands, the North and South islands. Finally, for $\hat\delta>\hat\delta_2$ we have three real roots, both in Eq.~\eqref{eq:cubic} and in Eq.~\ref{eq:cubic2}. On the $x$ axis, we can find the centers of three regions: The newly created West island, the stable center, and the usual East island, while we still have the North and South islands. Of course, for $\eps=0$ we retrieve, as $\hat\delta\to\infty$, the usual structure of the Hénon map.
\section{Use of the main tune to identify phase-space regions} \label{app:tune}
The frequency (main tune) associated to each closed orbit is evaluated using the so-called average phase advance method~\cite{Bartolini:292773,epac} over $N=4096$ turns, i.e.
\begin{equation}\nu_0 = \frac{1}{2\pi N}\sum_{n=1}^{N} \atan\frac{x_n p_{n-1} - x_{n-1} p_n}{x_n x_{n-1} + p_n p_{n-1}} \, .
\label{eq:tune1}
\end{equation}

\begin{figure}
    \includegraphics[width=0.5\linewidth,clip=]{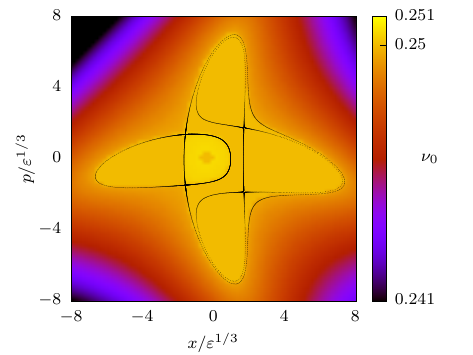}
    \includegraphics[width=0.5\linewidth,clip=]{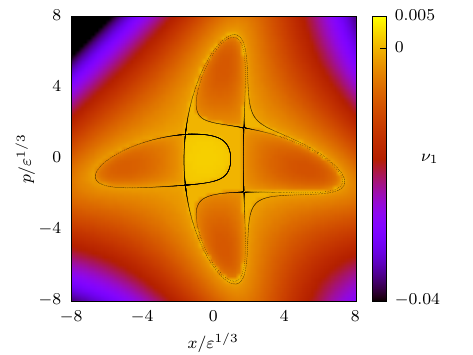}
    \caption{Left: Main tune $\nu_0$ (color scale, computed as in Eq.~\eqref{eq:tune1}) for a grid of initial conditions $(x_0,p_0)$ of the map of Eq.~\eqref{eq:doubleresmap}. Separatrices have also been represented (black line). Right: Secondary tune $\nu_1$ (color scale, computed as in Eq.~\eqref{eq:tune2}) for the same set of initial conditions and the same map (parameters values: $\Delta=0$, $\eps=10^{-4}$, $\delta/(2\pi)=0.005$, $\psi_0=0$).}
    \label{fig:tunemap}
\end{figure}

In the left plot of Fig.~\ref{fig:tunemap} we have computed the main tune $\nu_0$ for a set of initial conditions defined on a grid in the phase space $(x,p)$, and we observe that $\nu_0$ locks to the exact resonant value $\omegar/(2 \pi) =1/4$ not only within the islands, but also in the area enclosed by the outer separatrix. This result allows distinguishing between initial conditions that are in the region enclosed by these separatrices, but not those that are in the islands. To this end, we used the secondary tune $\nu_1$, which is the average phase advance computed on the stroboscopic map, \ie using only the $4$th iterate of the map, according to 
\begin{equation}\nu_1 = \frac{4}{2\pi N}\sum_{n=1}^{N/4} \atan\frac{x_{4n} p_{4n-4} - x_{4n-4}p_{4n}}{x_{4n} x_{4n-4} + p_{4n} p_{4n-4}}\, , 
\label{eq:tune2}
\end{equation}
and which is plotted, for the same set of initial conditions, in the right plot of Fig.~\ref{fig:phasespace_map}. In fact, it can be observed that $\nu_1$ changes sign within and outside the islands. The secondary tune $\nu_1$ is the rotation frequency of a particle around the center of the island and can be used to discriminate whether or not a particle is in an island. Moreover, its sign provides the direction of rotation around the stable fixed point, which can be used to distinguish between the islands and the core region.

Therefore, we define a particle to be trapped in an island if $\nu_0=1/4$ and $\nu_1 < 0$. Furthermore, by looking at the angle in the phase space of the final condition, it is possible to determine the specific island where the particle is. The particles with $\nu_1>0$ are in the core region and the remaining particles, \ie those with $\nu_0 \neq 1/4$ are classified in the external region. It should be mentioned that given the finite precision of the method used to determine $\nu_0$, a tolerance, based on the precision of the method~\cite{Bartolini:292773,epac}, should be defined to assess when $\nu_0$ can be considered equal to $1/4$.
\bibliography{biblio}
\bibliographystyle{unsrt}

\end{document}